\documentclass[12pt]{article}

\usepackage{amsmath}
\usepackage{amsfonts}
\usepackage{amsthm}
\usepackage[utf8]{inputenc}
\usepackage{lmodern}
\usepackage{titleps}
\usepackage{color}
\usepackage{slashbox}
\usepackage{multirow}

\date{}

\newtheorem{theo}{Theorem}
\newtheorem{lemma}{Lemma}

\newtheorem{prop}{Proposition}

\begin{document}

\title{Classical Sobolev orthogonal polynomials: eigenvalue problem}

\author{Juan F. Ma\~{n}as--Ma\~{n}as$^a$, Juan J. Moreno--Balc\'{a}zar$^{a,b}$.}

\maketitle

{\scriptsize
\noindent
$^a$Departamento de Matem\'{a}ticas, Universidad de Almer\'{\i}a, Spain.\\
$^b$Instituto Carlos I de F\'{\i}sica Te\'{o}rica y Computacional, Spain.\\

\noindent \textit{E-mail addresses:}  (jmm939@ual.es) J.F. Ma\~{n}as--Ma\~{n}as, (balcazar@ual.es) J.J. Moreno--Balc\'{a}zar.
}

\begin{abstract}

We consider the discrete Sobolev inner product
$$(f,g)_S=\int f(x)g(x)d\mu+Mf^{(j)}(c)g^{(j)}(c), \quad j\in \mathbb{N}\cup\{0\}, \quad c\in\mathbb{R}, \quad M>0, $$
where $\mu$ is a classical continuous measure with support on the real line (Jacobi, Laguerre or Hermite). The orthonormal polynomials with respect to this Sobolev inner product are eigenfunctions of a differential operator and  obtaining the asymptotic
behavior of the corresponding eigenvalues is the principal goal  of this paper.

\end{abstract}



 \section{Introduction}

The  classical continuous hypergeometric polynomials (CCHP) have been  well  known since the nineteenth century and they constitute a relevant class within the orthogonal polynomials.  Thus, almost all the books devoted to orthogonal polynomials and their applications included chapters or sections about CHHP, see among others \cite{ismail,sz}.  CCHP can be defined as the polynomial solutions to the hypergeometric differential equation \cite{Koekoek-book-hyper}
\begin{equation} \label{he}
\sigma(x)y''(x)+  \tau (x)y'(x)=\lambda_n y(x),
\end{equation}
where $\sigma $ and $\tau$ are polynomials with $\deg (\sigma)\le 2,$ $\deg (\tau) =1, $ and $\lambda_n=n\left(\tau'(x)+\frac{n-1}{2}\sigma''(x)\right)$. One can prove that these polynomial solutions  are orthogonal polynomials with respect to a weight function $w$. In fact, on the real line they are orthogonal with respect to a measure $\mu$ given by $d\mu(x)=w(x)dx$ where, up to affine transformations, $w$ corresponds to one of these situations: Jacobi case $w(x)=(1-x)^{\alpha}(1+x)^{\beta}$ with $\alpha, \beta>-1$ and $x\in (-1,1)$, Laguerre case $w(x)=x^{\alpha}e^{-x}$ with $\alpha>-1$ and $x\in (0,\infty)$, and Hermite case $w(x)=e^{-x^2}$ with $x$ on the real line.

The equation (\ref{he}) can be rewritten as $$\textbf{B}[y(x)]=\lambda_n y(x), $$
where $\textbf{B}$ is a differential operator defined as $\textbf{B}:=\sigma (x)\mathcal{D}^2+ \tau (x)\mathcal{D}, $ being $\mathcal{D}$ the usual derivative operator. In this way, the CCHP are the eigenfunctions of the operator $\textbf{B}$ and $\lambda_n$ are the corresponding eigenvalues. Both $\textbf{B}$ and $\lambda_n$ are explicitly known.

Since the second half of the last century an emergent theory of orthogonal polynomials in Sobolev spaces has risen, see for example the surveys \cite{Mar-moreno} and \cite{marxu}. The seminal papers on this topic linked Sobolev orthogonal polynomials (SOP) with the simultaneous approximation of a function and their derivatives, but now some recent applications have been found in \cite{ren-yu-2015,Takemura-et-2015}.

In this work we have considered a special case of SOP called discrete SOP which are orthogonal with respect to the Sobolev inner product
\begin{equation}\label{pro}
(f,g)_S=\int f(x)g(x)d\mu+Mf^{(j)}(c)g^{(j)}(c), \quad j\in \mathbb{N}\cup\{0\}, \quad M>0,
\end{equation}
where $\mu$ is a classical measure, i.e., $d\mu(x)=w(x)dx$ being $w$  one of the classical weights described previously. Notice that when $j=0$ we have the so-called Krall polynomials which were a first extension of the CCHP \cite{Krall-1980}.

The inner product (\ref{pro}) can be seen as a perturbation of the standard inner product $\int f(x)g(x)d\mu.$ Thus, it is natural to wonder how this perturbation influences on the corresponding orthogonal polynomials, for example, about the asymptotic behavior of these SOP. In fact,  there has been a wide literature about this so far (e.g. the previous surveys \cite{Mar-moreno} and \cite{marxu} and the references there in). On the one hand, it is obvious that the orthogonal polynomials with respect to (\ref{pro}) do not satisfy the hypergeometric equation (\ref{he}). However, in \cite{Jung-et-1997} the authors impose conditions so that the polynomials $q_n$ orthonormal with respect to (\ref{pro})   satisfy a (possibly infinite order)  differential equation.  Later, in \cite{Bavinck1998-85-95} and \cite{Bavinck1999} the differential operator is explicitly built as
$$
\mathbf{L}:=\sum_{i=1}^{\infty}r_i(x)\mathcal{D}^{i},
$$
where $r_i$  is a polynomial with $\deg (r_i)\leq i,$ and satisfying
$$
\mathbf{L}[q_n(x)]=\widetilde{\lambda}_n q_n(x).
$$

In this way, the orthonormal polynomials $q_n$ with respect to (\ref{pro}) are the eigenfunctions of the differential operator $\mathbf{L}$ and $\widetilde{\lambda}_n$ are the corresponding eigenvalues.

Thus,  we analyze the asymptotic behavior of the $\widetilde{\lambda}_n.$  We prove that this behavior is different of the one of $\lambda_n.$   A first approach of this problem was done in \cite{JAT-2018} although on that occasion the authors focused their attention on computing a value related to the convergence of a series in a left-definite space. Here, we tackle the problem in a wider framework.

 The structure of the paper is the following. In Section \ref{section2} we provide a brief background about eigenvalues of a differential operator related to  discrete Sobolev orthonormal polynomials. Sections \ref{section3} and \ref{section4} are devoted to obtaining the asymptotic behavior of the eigenvalues $\widetilde{\lambda}_n$ distinguishing two cases: symmetric and nonsymmetric. In both cases the technique is the same although, as we will see,  the computational details are different. Finally, in Section \ref{section5} we give  a summary table of the results and we comment them.

\section{Some known facts}\label{section2}
We consider a classical nonsymmetric  measure $\mu$ and the corresponding sequence of
orthonormal polynomials $\{p_n\}_{n\geq0}$, then it was established in \cite{Bavinck1998-85-95} and \cite{Bavinck1999} that when  $p_n^{(j)}(c)\neq0$ for all $n=j, j+1, j+2, \ldots$, we get \begin{equation}\label{def-autovalores-tilde}\widetilde{\lambda}_n=\lambda_n+M\alpha_n,
 \end{equation}
where $\lambda_n$ are the eigenvalues of the differential operator $\mathbf{B}$ and $\alpha_n$ is a sequence of real numbers such that if they are chosen conveniently, then the differential operator $\mathbf{L}$ is uniquely determined. To do this, it is enough to take   $\alpha_0=0$ and  $\{\alpha_i\}_{i=1}^{j}$  arbitrarily when $j>0.$

Thus, we have the measure $\mu$ such that $d\mu(x)=w(x)dx$ with $w$ corresponding to Laguerre o Jacobi case.  We take $c=0$ for the Laguerre case and $c\in \{-1,1\} $ for the Jacobi one.  Notice that the condition $p_n^{(j)}(c)\neq0$ is always satisfied when $c$ is chosen in this way because $\{p_n\}_{n\geq0}$ is a family of CCHP.
 Then, using \cite[f. (8-9)]{Bavinck1998-85-95} we get
\begin{equation}\label{alpha_n-caso-no-simetrico}
\alpha_n=\alpha_j+\sum_{i=j+1}^n(\lambda_i-\lambda_{i-1})K_{i-1}^{(j,j)}(c,c), \quad n\ge j+1,
\end{equation}
      where $K_{n}^{(j,j)}(c,c)$ denotes the partial derivatives of the $n$th kernel for the sequence of orthonormal polynomials $\{p_n\}_{n\ge 0}$ with respect to $\mu$, i.e.,

      $$K_n^{(r,s)}(x,y)=\sum_{i=0}^np_i^{(r)}(x)p_i^{(s)}(y),\qquad r,s\in\mathbb{N}\cup\{0\}.$$

Since  $\{\alpha_i\}_{i=1}^{j}$ can be chosen arbitrarily,  we take $\alpha_i=0$ for $i\in\{0,1,2,\dots,j\}$. Thus, (\ref{alpha_n-caso-no-simetrico}) is transformed into
\begin{equation}\label{alpha_n-caso-no-simetrico-2}
\alpha_n=\sum_{i=j+1}^n(\lambda_i-\lambda_{i-1})K_{i-1}^{(j,j)}(c,c), \quad  n\ge j+1.
\end{equation}

When $\mu$ is a classical symmetric measure with respect to the origin we take $c=0.$ According to \cite[Sect. 2.3]{Bavinck1999} to guarantee (\ref{def-autovalores-tilde}) it is necessary that $p_n^{(j)}(0)\neq0$ for all $n=j, j+1, j+2, \ldots$ with $n-j$ even.  But again this holds because $\{p_n\}_n$ is a family of CCHP.   Then,  applying the results in \cite[Sect. 2.3]{Bavinck1999} and taking into account again  that $\{\alpha_i\}_{i=1}^{j}$ are chosen arbitrarily,   we get
\begin{equation*}
\alpha_{j+2n}=\sum_{i=1}^n(\lambda_{j+2i}-\lambda_{j+2i-2})K_{j+2i-1}^{(j,j)}(0,0), \quad n\ge 1.
\end{equation*}

We remark that in this case, when $j$ is even the subsequence of orthonormal polynomials $\{q_{2n+1}\}_{n\ge 0}$ with respect to the discrete Sobolev inner product (\ref{pro}) matches  the one of standard orthonormal polynomials $\{p_{2n+1}\}_{n\ge 0}$. An analogous situation takes place when $j$ is odd, i.e.,  $\{q_{2n}\}_{n\ge 0}\equiv \{p_{2n}\}_{n\ge 0}.$

The Hermite case is an example of this situation, but we can also consider the Gegenbauer case which occurs when $\alpha=\beta$ in the Jacobi case. Thus, for this last case we can take $c$ as any value of the set $\{-1,0,1\}.$

With these results we are in condition to obtain the asymptotic behavior of the eigenvalues $\widetilde{\lambda}_n$ in the next section.

\section{Asymptotic behavior of the eigenvalues: the nonsymmetric case}\label{section3}

First, we give a joint approach to cases related to Jacobi and Laguerre weights. As we have mentioned in the previous section, we denote by  $\{p_n\}_{n\ge 0}$ the orthonormal polynomials with respect to the classical weights. We also use the  notation  $a_n\approx b_n$ meaning $\lim_{n\to+\infty}a_n/b_n=1.$ Then, it is easy to check that for these families of polynomials we have
\begin{equation}\label{pnk(c)}
p_n^{(k)}(c)\approx C_{k}(-1)^n n^{ak+b}, \quad 0\leq k\leq n, \quad  \mathrm{with} \quad 2(ak+b)+1>0,
\end{equation}
 where $C_k$ is a constant independent of $n.$ When we consider the Laguerre case then $c=0,$ and $c\in \{-1,1\}$ for the Jacobi case. For other nonclassical families satisfying (\ref{pnk(c)}) see \cite{dls2015}.

 It is worth noting that the factor $(-1)^n$ may appear or not, for example, it appears when $c=-1$ in the Jacobi case and it does not  when $c=1$ in the same case. However, for the results that we will obtain this factor will be not relevant from the asymptotic point of view.

\begin{lemma}\label{pro-com-asi-kernel}
Assuming the condition (\ref{pnk(c)}), we have
\begin{equation*}
\lim_{n\to+\infty}\frac{K_{n-1}^{(\ell,\ell)}(c,c)}{n^{2(a\ell+b)+1}}=\frac{C_{\ell}^2}{2(a\ell+b)+1}.
\end{equation*}
\end{lemma}
\noindent \textbf{Proof:} It is enough to use  Stolz's criterion and (\ref{pnk(c)}) to get
\begin{eqnarray*}
\lim_{n\to+\infty}\frac{K_{n-1}^{(\ell,\ell)}(c,c)}{n^{2(a\ell+b)+1}}
&=&\lim_{n\to+\infty}\frac{K_{n-1}^{(\ell,\ell)}(c,c)-K_{n-2}^{(\ell,\ell)}(c,c)}{n^{2(a\ell+b)+1}-(n-1)^{2(a\ell+b)+1}}\\
&=&\lim_{n\to+\infty}\frac{(p_{n-1}^{(\ell)}(c))^2}{(2(a\ell+b)+1)n^{2(a\ell+b)}}=\frac{C_{\ell}^2}{2(a\ell+b)+1}.  \qquad
\Box
\end{eqnarray*}

To obtain the asymptotic behavior of $\widetilde{\lambda}_n$ it is necessary to know the asymptotics of the sequence $\{\alpha_n\}$ and use (\ref{def-autovalores-tilde}). Thus, we establish it in the next result.

\begin{prop}\label{prop-1}
Assuming (\ref{pnk(c)}) and  $\lambda_n = \gamma n^2+\delta n$  with $\gamma,\delta\in\mathbb{R}$, then  we have

\begin{equation*}
  \begin{array}{ll}
    \displaystyle  \lim_{n\to+\infty}\frac{\alpha_n}{n^{2(aj+b)+3}}=\frac{2\gamma C_j^2}{(2(aj+b)+3)(2(aj+b)+1)}, & \hbox{if $\gamma\neq0$;} \\
&\\
    \displaystyle  \lim_{n\to+\infty}\frac{\alpha_n}{n^{2(aj+b+1)}}=\frac{\delta C_j^2}{2(aj+b+1)(2(aj+b)+1)}, & \hbox{if $\gamma=0$,}
  \end{array}
  \end{equation*}
where  $\alpha_n$ was defined in  (\ref{alpha_n-caso-no-simetrico-2}).
\end{prop}
\noindent \textbf{Proof:} Observe that
$$\lambda_n-\lambda_{n-1}=\left\{
                                              \begin{array}{ll}
                                                2\gamma n-\gamma+\delta, & \hbox{if $\gamma\neq0$;} \\
                                                \delta, & \hbox{if $\gamma=0$.}
                                              \end{array}
                                            \right.
$$
 Then, we need to distinguish two cases depending on $\gamma.$ In both situations we use Lemma \ref{pro-com-asi-kernel},  (\ref{alpha_n-caso-no-simetrico-2}), and again the Stolz's criterion to deduce the result:
\begin{itemize}
  \item If $\gamma\neq0$,
\begin{align*}
&\lim_{n\to+\infty}\frac{\alpha_n}{n^{2(aj+b)+3}}\\
=&\lim_{n\to+\infty}\frac{\sum_{i=j+1}^n(\lambda_i-\lambda_{i-1})K_{i-1}^{(j,j)}(c,c)-\sum_{i=j+1}^{n-1}(\lambda_i-\lambda_{i-1})K_{i-1}^{(j,j)}(c,c)}
{n^{2(aj+b)+3}-(n-1)^{2(aj+b)+3}}\\
=&\lim_{n\to+\infty}\frac{(\lambda_n-\lambda_{n-1})K_{n-1}^{(j,j)}(c,c)}
{(2(aj+b)+3)n^{2(aj+b+1)}}=\frac{2\gamma C_j^2}{(2(aj+b)+3)(2(aj+b)+1)}.
\end{align*}

  \item If $\gamma=0$,
\begin{eqnarray*}
\lim_{n\to+\infty}\frac{\alpha_n}{n^{2(aj+b+1)}}
&=&\lim_{n\to+\infty}\frac{(\lambda_n-\lambda_{n-1})K_{n-1}^{(j,j)}(c,c)}
{2(aj+b+1)n^{2(aj+b)+1}}\\
&=&\frac{\delta C_j^2}{2(aj+b+1)(2(aj+b)+1)}.\qquad \Box
\end{eqnarray*}
\end{itemize}

\begin{theo}\label{teo-case-1}
Let $\widetilde{\lambda}_n$ be the eigenvalues of the differential operator $\mathbf{L}$ related to the orthonormal polynomials $q_n$ with respect to (\ref{pro}). Under the hypothesis of Proposition \ref{prop-1}, we get

\begin{equation*}
  \begin{array}{ll}
    \displaystyle  \lim_{n\to+\infty}\frac{\widetilde{\lambda}_n}{n^{2(aj+b)+3}}=\frac{2\gamma MC_j^2}{(2(aj+b)+3)(2(aj+b)+1)}, & \hbox{if $\gamma\neq0$;} \\
&\\
    \displaystyle  \lim_{n\to+\infty}\frac{\widetilde{\lambda}_n}{n^{2(aj+b+1)}}=\frac{\delta MC_j^2}{2(aj+b+1)(2(aj+b)+1)}, & \hbox{if $\gamma=0$.}
  \end{array}
\end{equation*}
\end{theo}

\noindent \textbf{Proof:} We only need to take limits in (\ref{def-autovalores-tilde}) and apply Proposition \ref{prop-1}. We only show the proof when  $\gamma\neq0$, the other case is totally similar.

\begin{eqnarray*}
\lim_{n\to+\infty}\frac{\widetilde{\lambda}_n}{n^{2(aj+b)+3}}
&=&\lim_{n\to+\infty}\frac{\lambda_n+M\alpha_n}{n^{2(aj+b)+3}}\\
&=&\lim_{n\to+\infty}\frac{\lambda_n}{n^{2(aj+b)+3}}+M\lim_{n\to+\infty}\frac{\alpha_n}{n^{2(aj+b)+3}}\\
&=&\frac{2\gamma M C_j^2}{(2(aj+b)+3)(2(aj+b)+1)}.
\end{eqnarray*}
The first limit is 0 because using (\ref{pnk(c)})  we have   $2(aj+b)+3>2.$ $\qquad \Box$

\subsection{Discrete Jacobi--Sobolev case}

We consider the discrete Sobolev inner product

\begin{equation}\label{pro-jacobi}
(f,g)_{JS}=\int_{-1}^1f(x)g(x)(1-x)^{\alpha}(1+x)^{\beta}dx+Mf^{(j)}(1)g^{(j)}(1),
\end{equation}
with $\alpha, \beta>-1$ and $j\in\mathbb{N}\cup\{0\}.$ We denote by $\{p_n^{(\alpha,\beta)}\}_{n\ge 0}$ the sequence of the classical Jacobi orthonormal polynomials with respect to the weight function $(1-x)^{\alpha}(1+x)^{\beta}.$  This inner product corresponds to (\ref{pro}) with $c=1.$

Using the properties of Jacobi polynomials (e.g., see \cite[f. (4.1.1), (4.3.3), (4.21.7)]{sz}), we deduce
\begin{equation*}
\left(p_n^{(\alpha,\beta)}\right)^{(k)}(1)\approx \frac{1}{2^{k+\frac{\alpha+\beta}{2}}\Gamma({\alpha+k+1})}n^{2k+\alpha+1/2},
\end{equation*}
so,  (\ref{pnk(c)}) is satisfied with $C_k=\frac{1}{2^{k+\frac{\alpha+\beta}{2}}\Gamma({\alpha+k+1})}$, $a=2$ and $b=\alpha+1/2.$   Since $\alpha>-1$, the condition $2(ak+b)+1>0$  holds.

On the other hand, Jacobi polynomials $p_n^{(\alpha,\beta)}$ satisfy the second--order differential equation (e.g., see \cite[f. (4.2.1)]{sz}):
$$(x^2-1)y''(x)+(\alpha-\beta+(\alpha+\beta+2)x)y'(x)=n(n+\alpha+\beta+1)y(x),$$
 thus, we deduce $\lambda_n=n^2+n(\alpha+\beta+1)$.

Now, we are ready to apply Theorem \ref{teo-case-1},  getting
$$\lim_{n\to+\infty}\frac{\widetilde{\lambda}_n}{n^{4j+2\alpha+4}}
=\frac{M}{2^{2j+\alpha+\beta+1}(2j+\alpha+2)(2j+\alpha+1)\Gamma^2(\alpha+j+1)}.$$
A similar result can be obtained if we choose $c=-1$ in (\ref{pro-jacobi}) instead of $c=1.$

\subsection{Discrete Laguerre--Sobolev case}

Now, we consider
\begin{equation*}
(f,g)_{LS}=\int_{0}^{+\infty}f(x)g(x)x^{\alpha}e^{-x}dx+Mf^{(j)}(0)g^{(j)}(0), \quad \alpha>-1, \quad j\in\mathbb{N}\cup\{0\}.
\end{equation*}
  We denote by $\{l_n^{(\alpha)}\}_{n\ge 0}$ the sequence of the classical Laguerre orthonormal polynomials with respect to the weight function $x^{\alpha}e^{-x}.$  We have taken $c=0$ in (\ref{pro}).
In this case we use the properties of Laguerre polynomials (e.g., see \cite[f. (5.1.1), (5.1.7), (5.1.14)]{sz}) to obtain
\begin{equation*}
\left(l_n^{(\alpha)}\right)^{(k)}(0)\approx\frac{(-1)^k}{\Gamma(\alpha+k+1)}n^{k+\alpha/2}.
\end{equation*}
Again, (\ref{pnk(c)}) is satisfied taking now $C_k=\frac{(-1)^k}{\Gamma(\alpha+k+1)}$, $a=1$ and $b=\alpha/2.$  Since classical Laguerre  polynomials satisfy the hypergeometric differential equation (e.g., see \cite[f. (5.1.2)]{sz}):
$$-xy''(x)+(x-\alpha-1)y'(x)=ny(x),$$
we have $\lambda_n=n$. Therefore, we can apply Theorem \ref{teo-case-1} taking into account that in this case $\gamma=0$ and $\delta=1$, getting
 $$\lim_{n\to+\infty}\frac{\widetilde{\lambda}_n}{n^{2j+\alpha+2}}=\frac{M}{(2j+\alpha+2)(2j+\alpha+1)\Gamma^2(\alpha+j+1)}.$$

\section{Asymptotic behavior of the eigenvalues: the symmetric case}\label{section4}

We suppose that  $\mu$ is a symmetric measure and we take $c=0$. Thus we can proceed like in the previous section, but bearing in mind that now both families of orthonormal polynomials, $\{p_n\}$ and $\{q_n\}$, are symmetric.
Therefore, we have to assume  similar conditions to (\ref{pnk(c)}) for the subsequences of even and odd polynomials. Thus, we suppose that
\begin{eqnarray}\label{pnk(c)-simetrico1}
p_{2n}^{(2k)}(0)\approx C_{k,1}(-1)^n n^{a_1k+b_1}, \quad  p_{2n+1}^{(2k+1)}(0)\approx C_{k,2}(-1)^n n^{a_2k+b_2},
\end{eqnarray}
with $2(a_1k+b_1)+1>0$ and $2(a_2k+b_2)+1>0$ for all $k\in \{0, \ldots,  n\}.$

Assuming (\ref{pnk(c)-simetrico1}), we can obtain similar results to the ones obtained in the previous section. Since the techniques are the same we only state the main outcome.
\begin{theo}\label{teo-case-2}
Let $\widetilde{\lambda}_n$ be the eigenvalues of the differential operator $\mathbf{L}$ related to the orthonormal polynomials $q_n$ with respect to  (\ref{pro}). We assume (\ref{pnk(c)-simetrico1}) and $\lambda_n= \gamma n^2+\delta n$  with $\gamma,\delta\in\mathbb{R}.$ Then,
\begin{itemize}
\item  If $j=2r$, we get
\begin{equation*}
  \begin{array}{ll}
    \displaystyle  \lim_{n\to+\infty}\frac{\widetilde{\lambda}_{2r+2n}}{n^{2(a_1r+b_1)+3}}
=\frac{8\gamma MC_{r,1}^2}{(2(a_1r+b_1)+3)(2(a_1r+b_1)+1)}, & \hbox{if $\gamma\neq0$;} \\
&\\
    \displaystyle  \lim_{n\to+\infty}\frac{\widetilde{\lambda}_{2r+2n}}{n^{2(a_1r+b_1+1)}}=
\frac{\delta MC_{r,1}^2}{(a_1r+b_1+1)(2(a_1r+b_1)+1)}, & \hbox{if $\gamma=0$.}
  \end{array}
\end{equation*}
\item If $j=2r+1$, we get
\begin{equation*}
  \begin{array}{ll}
    \displaystyle  \lim_{n\to+\infty}\frac{\widetilde{\lambda}_{2r+1+2n}}{n^{2(a_2r+b_2)+3}}
=\frac{8\gamma MC_{r,2}^2}{(2(a_2r+b_2)+3)(2(a_2r+b_2)+1)}, & \hbox{if $\gamma\neq0$;} \\
&\\
    \displaystyle  \lim_{n\to+\infty}\frac{\widetilde{\lambda}_{2r+1+2n}}{n^{2(a_2r+b_2+1)}}=
\frac{\delta MC_{r,2}^2}{(a_2r+b_2+1)(2(a_2r+b_2)+1)}, & \hbox{if $\gamma=0$.}
  \end{array}
\end{equation*}
\end{itemize}
\end{theo}

\subsection{Discrete Hermite--Sobolev case}

We take the Hermite weight function and $c=0$, then the inner product (\ref{pro}) is transformed into
\begin{equation*}
(f,g)_{HS}=\int_{-\infty}^{+\infty}f(x)g(x)e^{-x^2}dx+Mf^{(j)}(0)g^{(j)}(0), \quad j\in\mathbb{N}\cup\{0\}.
\end{equation*}
 We denote by $\{h_n\}_{n\ge 0}$ the sequence of the classical Hermite orthonormal polynomials with respect to the weight function $e^{-x^2}.$ Using the properties of Hermite polynomials
(e.g., see \cite[f. (5.5.1), (5.5.5), (5.5.10)]{sz}) we deduce
\begin{eqnarray*}
h_{2n}^{(2k)}(0)\approx (-1)^n\frac{(-1)^{k}2^{2k}}{\sqrt{\pi}}n^{k-1/4},\quad
h_{2n+1}^{(2k+1)}(0)\approx (-1)^n\frac{(-1)^{k}2^{2k+1}}{\sqrt{\pi}}n^{k+1/4}.
\end{eqnarray*}
  Therefore, (\ref{pnk(c)-simetrico1}) holds with $C_{k,1}=\frac{(-1)^{k}2^{2k}}{\sqrt{\pi}}$, $a_1=1$, $b_1=-1/4$, $C_{k,2}=\frac{(-1)^{k}2^{2k+1}}{\sqrt{\pi}}$, $a_2=1$ and $b_2=1/4$.

Moreover, Hermite polynomials satisfy the second--order differential equation (e.g., see \cite[f. (5.5.2)]{sz})
$$-y''(x)+2xy'(x)=2ny(x).$$
Then, we have $\lambda_n=2n$. In this way, we can apply  Theorem \ref{teo-case-2}  with $\gamma=0$ and $\delta=2$, obtaining the corresponding asymptotic behavior of the eigenvalues $\widetilde{\lambda}_n$, i.e.,

\begin{itemize}
\item If $j=2r$, then
$$\lim_{n\to+\infty}\frac{\widetilde{\lambda}_{2r+2n}}{n^{2r+3/2}}=\frac{M2^{4r+1}}{\pi(r+3/4)(2r+1/2)}.$$
\item If $j=2r+1$, then
$$\lim_{n\to+\infty}\frac{\widetilde{\lambda}_{2r+1+2n}}{n^{2r+5/2}}=\frac{M2^{4r+3}}{\pi(r+5/4)(2r+3/2)}.$$
\end{itemize}

\subsection{Discrete Gegenbauer--Sobolev case}
For this case, we consider the discrete Sobolev inner product
\begin{equation*}
(f,g)_{GS}=\int_{-1}^{1}f(x)g(x)(1-x^2)^{\alpha}dx+Mf^{(j)}(0)g^{(j)}(0), \quad  \alpha>-1, \quad j\in\mathbb{N}\cup\{0\}.
\end{equation*}
We denote by $\{c_n^{(\alpha)}\}_{n\ge 0}$ the sequence of the classical Gegenbauer orthonormal polynomials with respect to the weight function $(1-x^2)^{\alpha}.$  Using some properties of these polynomials (e.g.,  \cite[f. (4.7.1), (4.7.14), (4.7.15), (4.7.30)]{sz} and the relations in \cite[p. 60]{sz} for $\alpha=-1/2$; notice that in \cite{sz} the author works with $c_n^{(\lambda-1/2)}$ with $\lambda>-1/2$), we get
\begin{eqnarray*}
\left(c_{2n}^{(\alpha)}\right)^{(2k)}(0)&\approx&(-1)^{n}\frac{(-1)^k 2^{2k+1/2}}{\sqrt{\pi}}n^{2k},\\
\left(c_{2n+1}^{(\alpha)}\right)^{(2k+1)}(0)&\approx&(-1)^{n}\frac{(-1)^k 2^{2k+3/2}}{\sqrt{\pi}}n^{2k+1}.
\end{eqnarray*}
Then,  (\ref{pnk(c)-simetrico1}) holds with $C_{k,1}=\frac{(-1)^k 2^{2k+1/2}}{\sqrt{\pi}},$
$a_1=2,$
$b_1=0,$
$C_{k,2}=\frac{(-1)^k 2^{2k+3/2}}{\sqrt{\pi}},$
$a_2=2,$
and $b_2=1. $ This family of polynomials satisfies the hypergeometric differential equation (e.g., see \cite[f. (9.8.23)]{Koekoek-book-hyper})
$$(x^2-1)\left(c_n^{(\alpha)}\right)''(x)+2(\alpha+1)x\left(c_n^{(\alpha)}\right)'(x)=n(n+2\alpha+1)c_n^{(\alpha)}(x),$$
 so, we have $\lambda_n=n^2+n(2\alpha+1)$. Finally, we apply Theorem \ref{teo-case-2} with  $\gamma=1$ and $\delta=2\alpha+1$, getting

\begin{itemize}
\item If $j=2r$,
$$\lim_{n\to+\infty}\frac{\widetilde{\lambda}_{2r+2n}}{n^{4r+3}}
=\frac{2^{4(r+1)}M}{\pi(4r+3)(4r+1)},$$
\item If $j=2r+1$,
$$\lim_{n\to+\infty}\frac{\widetilde{\lambda}_{2r+1+2n}}{n^{4r+5}}
=\frac{2^{4r+6}M}{\pi(4r+5)(4r+3)}.$$
\end{itemize}

\section{Conclusions}\label{section5}

Finally, we provide  a summary table of the results obtained for SOP and we compare them with those ones known for classical polynomials.
\begin{center}
\begin{tabular}{|c|*{2}{c|}}\hline
\backslashbox{Eigenvalues}{Case}&Jacobi&Laguerre\\\hline
Asymptotics of $\lambda_n$              & $n^2$                            & $n$                             \\  \hline
Asymptotics of $\widetilde{\lambda}_n$  & $\mathcal{C}_1n^{4j+2\alpha+4}$  & $\mathcal{C}_2n^{2j+\alpha+2}$  \\\hline
\end{tabular}
\end{center}

\begin{center}
\begin{tabular}{|c|c|c|c|} \cline{2-4}
\multicolumn{1}{c|}{ } & \backslashbox{Eigenvalues}{Case} & Hermite & Gegenbauer \\ \cline{2-4}
\multicolumn{1}{c|}{ } & Asymptotics of $\lambda_n$ & $2n$ & $n^2$  \\  \hline
\multirow{2}{2.5cm}{ If $j=2r$ }
& Asymptotics of $\widetilde{\lambda}_{2n}$         & $\mathcal{C}_3 n^{2r+3/2}$ & $\mathcal{C}_4 n^{4r+3}$    \\ \cline{2-4}
& Asymptotics of $\widetilde{\lambda}_{2n+1}$        & $4n$ & $4n^2$    \\ \hline
\multirow{2}{2.5cm}{ If $j=2r+1$ }
& Asymptotics of $\widetilde{\lambda}_{2n}$                 & $4n$ & $4n^2$    \\ \cline{2-4}
& Asymptotics of $\widetilde{\lambda}_{2n+1}$    & $\mathcal{C}_5 n^{2r+5/2}$ & $\mathcal{C}_6 n^{4r+5}$     \\ \hline
\end{tabular}
\end{center}
The constants $\mathcal{C}_i,\, i=1, \ldots, 6.$  can be found explicitly in the previous sections.
\medskip

 We can observe that in all the cases the presence of the discrete part $Mf^{(j)}(c)g^{(j)}(c)$ in the Sobolev inner product  leads to important changes in the asymptotic behavior of the eigenvalues. For example, in the Laguerre case the  growing orden  of the eigenvalues increases $2j+\alpha+1>0$, i.e., when $M=0$ the eigenvalues have a linear growth that  changes to a growth of $O(n^{2j+\alpha+2})$ if $M>0.$ Similar situations occur in the rest of the cases. We can observe that this change is bigger in the bounded cases: Jacobi and Gegenbauer.

\bigskip

\noindent
\textbf{Acknowledgments.} The authors are partially supported by the Ministry of Science, Innovation and Universities of Spain and the European Regional Development Fund (ERDF), grant MTM2017-89941-P and by Research Group FQM-0229 (belonging to Campus of International Excellence CEIMAR). The author JFMM is funded by a grant of Plan Propio de la Universidad de Almer\'{\i}a. The author JJMB is also partially supported by the research centre CDTIME of Universidad de Almer\'{\i}a and  by Junta de Andaluc\'{\i}a and ERDF, ref. SOMM17/6105/UGR.

\end{document}